\documentclass{amsart}
\usepackage{amsmath,amssymb,amsthm,url,color}
\newtheorem{thr}{Theorem}
\newtheorem{lem}[thr]{Lemma}

\theoremstyle{definition}

\newtheorem{prob}[thr]{Problem}

\theoremstyle{remark}

\newtheorem{obs}[thr]{Observation}

\def\F{\mathcal{F}}
\def\R{\mathbb{R}}
\def\A{\mathcal{A}}

\begin{document}

\title{Nonnegative rank depends on the field}

\author{Yaroslav Shitov}

\address{National Research University Higher School of Economics, 20 Myasnitskaya ulitsa, Moscow 101000, Russia}
\email{yaroslav-shitov@yandex.ru}


\begin{abstract}
We present an example of a subfield $\F\subset\R$ and a matrix $A$ whose conventional and nonnegative ranks equal five, but the nonnegative rank with respect to $\F$ equals six. In other words, $A$ can be represented as a sum of five rank-one matrices with nonnegative real entries but not as a sum of five rank-one matrices with nonnegative entries in $\F$.
\end{abstract}

\maketitle

\section{Introduction}

Let $A$ be an $m\times n$ nonnegative matrix, that is, a matrix with nonnegative real entries. The \textit{nonnegative rank} of $A$ is the smallest $k$ for which $A$ can be written as a product $BC$ of two nonnegative matrices of sizes $m\times k$ and $k\times n$, respectively. This concept is a part of the \textit{nonnegative matrix factorization} problem (NMF), which is essentially the task to find an optimal approximation of $A$ with a product as above. This question is important in those topics of modern science that naturally deal with nonnegative arrays of data, because in such cases the output of NMF can be easier to interpret in comparison with other factorization techniques~\cite{Gillis}, giving the applications in the areas of quantum mechanics~\cite{CR}, image processing~\cite{LS}, statistics~\cite{KRS}, text mining~\cite{textmin}, music analysis~\cite{musan}, and many other topics.

Another striking application of NMF comes from the geometric point of view and lies in the optimization perspective. The smallest number of inequalities needed to define a given polytope $P$ up to a linear projection can be expressed as the nonnegative rank of the so-called \textit{slack matrix} of $P$, so the small-rank nonnegative factorizations of this matrix may lead to compact linear programming formulations and allow a faster optimization over $P$, see~\cite{FMPTdW, Yan}. Also, the geometric perspective gives a better understanding of the complexity of algorithms solving the NMF problem~\cite{Moitra, mypsd} for upper bounds on the complexity and~\cite{Vav} for the NP-hardness proof by Vavasis.

Although a complete understanding of the algorithmic complexity of NMF in its general formulation is now reached (the problem is $\exists\R$-complete, that is, equivalent to deciding if a given system of polynomial equations with integral coefficients has a real solution~\cite{myuni}), many important questions remain open. Is there a polynomial-time algorithm computing the nonnegative rank of a matrix with bounded rank~\cite{GG}? Another open issue is the \textsc{exact nmf} problem, which asks whether or not the nonnegative rank of a given matrix equals its conventional rank. Vavasis~\cite{Vav} proved the NP-hardness of this problem and asked whether or not it belongs to NP. The answer to this question would be 'yes' if the property of having nonnegative rank less than a given integer could be certified by a factorization whose entries are polynomial-time computable rational functions of the entries of the initial matrix. In our paper, we show the non-existence of such factorizations in general, but since the certificates involved in potential NP descriptions of \textsc{exact nmf} are not \textit{a priori} limited to explicit factorizations of input matrices, our result does not rule out the possibility that \textsc{exact nmf} is NP-complete.

More precisely, we prove the following result. We provide an example of a subfield $\F\subset\R$ and a matrix $A$ with nonnegative entries in $\F$ such that
\begin{equation}\label{eq5}
5=\mathrm{rank}(A)=\mathrm{rank}_+(A,\R)<\mathrm{rank}_+(A,\F)=6,
\end{equation}
where $\mathrm{rank}_+(A,\mathcal{K})$ is the smallest integer $k$ such that $A$ can be written as a sum of $k$ rank-one matrices with nonnegative entries in a field $\mathcal{K}$. This value is called the nonnegative rank of $A$ \textit{with respect to} $\mathcal{K}$, and, in particular, the function $\mathrm{rank}_+(A,\R)$ corresponds to the usual nonnegative rank as defined above. Translating our result into the geometric language, we get the following.

\begin{thr}\label{prPolytope}
There exist a field $\F\subset\R$ and polytopes $P,Q\subset\mathbb{R}^{4}$ such that:

\noindent (I) the vertices of both $P$ and $Q$ have coordinates in $\mathcal{F}$;

\noindent (II) $\dim P=4$ and $P\subset Q$;

\noindent (III) there is a $4$-simplex $\Delta$ satisfying $P\subset\Delta\subset Q$,

\noindent (IV) every such $\Delta$ has a vertex not all of whose coordinates belong to $\mathcal{F}$.
\end{thr}

Our construction appeared in~\cite{mydep}, and it was the earliest known solution of the problem by Cohen and Rothblum, who asked in 1993 whether the nonnegative rank of a matrix can be sensitive to the ordered field with respect to which it is being computed. Despite the subsequent advances, our construction remains an only known proof that the answer to \textsc{exact nmf} may depend on the underlying field, --- in particular, the papers~\cite{anotherproof, myuni, myCR} allow to distinguish between the nonnegative ranks with respect to the more conventional setting of the fields $\mathbb{Q}$ and $\mathbb{R}$ but do not consider the case $\mathrm{rank}(A)=\mathrm{rank}_+(A,\R)$ at all, while the studies~\cite{KK, Consp} work on $\mathrm{rank}(A)=\mathrm{rank}_+(A,\R)=4$ without definitive progress on the problem being discussed.

\medskip

In the rest of this note, we follow the geometrical approach and prove Theorem~\ref{prPolytope}, so let us comment on why does this theorem imply the existence of a matrix $A$ as in~\eqref{eq5}. In fact, such a matrix can be constructed as what is called the \textit{slack matrix} of the pair $(P,Q)$, that is, the matrix which has rows indexed by vertices of $P$, columns indexed by facets of $Q$, and an $(i,j)$ entry equal to the distance from the $i$th vertex to the $j$th facet. The details of the proof of this fact are available in the first arXiv version~\cite{mydep} of this paper, but since this result can be proved by standard techniques, we do not reproduce the proof here. Also, the desired implication can be seen as a consequence of the equality between the \textit{extension complexity} of a polytope and the nonnegative rank of its slack matrix~\cite{FMPTdW, Yan}. Alternatively, the fact that Theorem~\ref{prPolytope} implies~\eqref{eq5} can be deduced from the result of Vavasis~\cite{Vav} on the equivalence of \textsc{exact nmf} and \textsc{intermediate simplex}. Let us also mention that the combinatorial approach to bounding the nonnegative rank may not be of help to obtain similar results. In particular, the vertices of the polytope $P$ in our construction lie in the interior of $Q$, so the corresponding slack matrix has all entries positive, which makes several combinatorial lower bounds for nonnegative rank as in~\cite{Fior} inapplicable to our example. This contrasts the situation in the paper~\cite{myCR} making a heavy use of a zero-nonzero pattern of a matrix to derive lower bounds on its nonnegative rank.

Our paper is organized as follows. In Section~2, we construct the field $\F$ as in Theorem~\ref{prPolytope} and compute several numbers used in our description of the polytopes $P,Q,\Delta$. In Section~3, we give a formal description of these polytopes, and we check the items (I)--(III) in Theorem~\ref{prPolytope}. In Section~4, we complete our proof by checking the item (IV) in Theorem~\ref{prPolytope}. Our paper requires some computer calculations, and we will sometimes refer to our \textit{Wolfram Mathematica} files~\cite{files} containing all necessary numerical data.

\section{The field $\F$}\label{seccompdet}

In this section, we explain how to compute the (unique) numbers $\alpha=-0.0311...$, $\beta=-0.4088...$, $\gamma=0.3983...$ that make the polynomial
$\pi(t)=(-96 \alpha - 96 \beta - 48 \alpha \beta - 192 \gamma - 96 \alpha \gamma -
  96 \beta \gamma + 99 \alpha \beta \gamma)+ $$ $$t(- 1392 \alpha - 1296 \beta -
  96 \alpha \beta - 1696 \gamma - 1624 \alpha \gamma -
  952 \beta \gamma + 1370 \alpha \beta \gamma ) +$$ $$t^2( - 256  -
  8608 \alpha - 5488 \beta + 1172 \alpha \beta -
  5536 \gamma - 10240 \alpha \gamma - 2164 \beta \gamma +
  4292 \alpha \beta \gamma )+$$ $$t^3( - 2944  - 25648 \alpha -
  8512 \beta  + 4100 \alpha \beta - 8192 \gamma -
  21056 \alpha \gamma - 1736 \beta \gamma +
  4768 \alpha \beta \gamma ) +$$ $$t^4 ( - 8576 - 32288 \alpha -
  4832 \beta + 4072 \alpha \beta - 2912 \gamma  -
  12608 \alpha \gamma - 440 \beta \gamma +
  1744 \alpha \beta \gamma ) +$$ $$t^5( - 3584  - 11648 \alpha -
  896 \beta + 1120 \alpha \beta )\in\mathbb{Q}[t]$
have the form $(-t-2)u^2$, where $u\in\mathbb{R}[t]$ has degree two. Our ground field $\mathcal{F}$ is taken as $\mathbb{Q}(\alpha,\beta,\gamma)$, which is the smallest subfield of $\mathbb{R}$ containing $\alpha,\beta,\gamma$. We define $\tau=0.1765...$ to be one of the two roots of $u$, and later will need to check that $\tau\notin\F$.

A straightforward approach to computing $\alpha,\beta,\gamma$ is based on the fact that they are a solution of the polynomial system
$$\pi(-2)=0,\,\,s_0(\pi,\partial\pi/\partial t)=0,\,\,s_1(\pi,\partial\pi/\partial t)=0,$$
where $s_i(f,g)$ denotes the $i$th subresultant of polynomials $f,g$. Unfortunately, we did not manage to make this method work because of its algorithmic complexity, so we will use a different method with all the details provided in File~1 of~\cite{files}. Namely, we define a substitution $t_0=-t-2$ and set $\pi_0(t_0)=\pi(-2-t_0)$, and we note that $\pi$ can be written in the desired form $(-t-2)u^2$ if and only if $\pi_0(0)=0$ and $\pi_0/t_0$ is a square of a quadratic polynomial. The condition $\pi_0(0)=0$ allows us to express $\gamma$ as a rational function of $\alpha,\beta$, so we can replace every appearance of $\gamma$ in $\pi_0/t_0$ by this function. We deal with the condition of $\pi_0/t_0$ being a square as follows.

\begin{lem}\label{lemsq}
The polynomial $y_4 t^4  +  y_3 t^3  +  y_2 t^2  +  y_1 t  +  y_0$ is a square of a quadratic polynomial if and only if
$$\frac{y_3^3}{8 y_4^3} - \frac{y_2 y_3}{2 y_4^2} + \frac{y_1}{y_4} =0,\,\,\,
\left( \frac{y_2 y_3^2}{4y_4^3}-\frac{3 y_3^4}{64y_4^4} - \frac{y_1 y_3}{y_4^2} + \frac{4y_0}{y_4}\right) - \left(\frac{3 y_3^2}{8 y_4^2} - \frac{y_2}{y_4}\right)^2 = 0.$$
In this case, its roots are
$$\pm\sqrt{\frac{3 y_3^2}{16 y_4^2}-\frac{y_2}{2 y_4}} - \frac{y_3}{4 y_4}.$$
\end{lem}

\begin{proof}
See our computation, based on the idea from~\cite{Spey}, in File~2 of~\cite{files}.
\end{proof}

We apply this lemma to the polynomial $\pi_0/t_0$, where the $y_i$'s are the corresponding rational functions depending on $\alpha$ and $\beta$. This allows us to compute $\alpha, \beta$, and then we can find $\gamma, \tau$. Finally, we use standard tools to check that $\F=\mathbb{Q}(\alpha,\beta,\gamma)$ has degree $16$ over $\mathbb{Q}$ while  $\mathbb{Q}(\alpha,\beta,\gamma,\tau)$ has degree $32$ over $\mathbb{Q}$, so we get $\tau\notin\F$.

\section{The polytopes $P$, $Q$, $\Delta$}

This section gives a formal description of the polytopes $P$, $Q$, $\Delta$. Our construction is four-dimensional, and we work in the affine sub\-space $\A=\{x_1+\ldots+x_5=1\}$ in $\mathbb{R}^5$; we say that four points in $\A$ are in \textit{general position} if their affine span is a \textit{hyperplane}, that is, a three-dimensional plane in $\A$.

\sloppy{We set $\Omega(t)=\mu_t(1 + t, 1 + 2 t,1 + t,1,0)$, $F_{11}=\lambda_{11}(1000 + 514 \alpha, 1056, 524 + 131 \alpha, 772 + 193 \alpha, 648 +
 162 \alpha)$, $F_{12}=\lambda_{12}(1000 + 532 \alpha, 1128, 1012 + 253 \alpha, 236 + 59 \alpha, 624 +
 156 \alpha)$, $F_{13}=\lambda_{13}(500 + 233 \alpha, 432, 536 + 134 \alpha, 176 + 44 \alpha, 356 + 89 \alpha)$,
$F_{21}=\lambda_{21}(14044 + 3511 \beta, 20000 + 9467 \beta, 17868, 8888 +2222 \beta, 19200 + 4800 \beta)$,
$F_{22}=\lambda_{22}(136 + 34 \beta, 200 + 98 \beta, 192, 144 + 36 \beta, 128 + 32 \beta)$,
$F_{23}=\lambda_{23}(96 + 24 \beta, 100 + 28 \beta, 12, 72 + 18 \beta, 120 + 30 \beta)$,
$F_{31}=\lambda_{31}(26 - 3 \gamma, 22 - 2 \gamma, 25 + 5 \gamma, 20, 7)$,
$F_{32}=\lambda_{32}(17 - 3 \gamma, 22 - 2 \gamma, 25 + 5 \gamma, 20, 16)$,
$F_{33}=\lambda_{33}(376 - 81 \gamma, 384 - 54 \gamma, 500 + 135 \gamma, 540, 200)$,
$F_{41}=\lambda_{41}(618, 392, 365, 625, 500)$,
$F_{42}=\lambda_{42}(1863, 1252, 1250, 1875, 1260)$,
$F_{43}=\lambda_{43}(384, 496, 495, 625, 500)$, $H=(3,3,3,3,4)/16$,
where the coefficients $\mu_t$ and $\lambda$ are such that the corresponding vectors are \textit{normalized}, that is, their coordinates sum to $1$.}

In the rest of our paper, $\varepsilon$ denotes a sufficiently small positive number. We consider a rational approximation of $\tau$, that is, a pair $(q_1,q_2)$ of rational numbers satisfying $q_1<\tau<q_2$ and $|q_1-q_2|<\varepsilon$. We define $\Omega$ as the point $\Omega(\tau)$ obtained by plugging in $t=\tau$ in the function $\Omega(t)$ defined above, and, similarly, we set $\Omega_i=\Omega(q_i)$. We define $\omega$ as a generic point in $\A$ (and we think of it as a point close to $\Omega$), and we consider the functions $f_{ij}:\R\to\A$ defined by $f_{ij}=(F_{ij}-(0,0,0,0,v_{ij}))/(1-v_{ij})$, where the $v_{ij}$'s are variables. In other words, we consider the straight line $f_{ij}(\mathbb{R})$, which comes as the intersection of $\mathcal{A}$ and the $2$-plane containing the origin and spanned by the vectors $F_{ij}$ and $(0,0,0,0,1)$. Also, we denote the list of all $v_{ij}$'s by $v=(v_{ij})$.

For every $i\in\{1,2,3,4\}$, we check (see File~3 in~\cite{files}) that the points $\Omega,F_{i1},F_{i2},F_{i3}$ are in general position, so the points $\omega, f_{i1}, f_{i2}, f_{i3}$ are still in general position whenever $(\omega,v)$ is sufficiently close to $(\Omega,0,\ldots,0)$. In this case, there is, up to scaling, a unique vector $\pi_i(\omega,v)$ such that $\pi_i(\omega,v)\cdot x=0$ is the equation defining the hyperplane passing through $\omega, f_{i1}, f_{i2}, f_{i3}$. We define $V_j(\omega,v)$ as the point that lies at the intersection of the hyperplane $\{x_j=0\}$ and all the hyperplanes $\{\pi_i(\omega,v)\cdot x=0\}$ with $i$ different from $j$. We define $\pi_5(\omega,v)$ as a vector satisfying $\pi_5(\omega,v)\cdot x=0$ whenever $x$ belongs to the affine hull of the points $V_1(\omega,v),V_2(\omega,v),V_3(\omega,v),V_4(\omega,v)$. The convex hull of these points together with $\omega$ is denoted by $\Delta(\omega,v)$. We write $\pi_k=\pi_k(\Omega,0,\ldots,0)$, $V_j=V_j(\Omega,0,\ldots,0)$, and $\Delta=\Delta(\Omega,0,\ldots,0)$.

\begin{lem}\label{lempi14}
The points $\Omega,V_1,V_2,V_3,V_4$ are affinely independent, so $\Delta$ is indeed a simplex. We can choose $\pi_i$ to be normalized, and then $\pi_i\cdot x\geqslant 0$ are the inequalities defining $\Delta$.
\end{lem}

\begin{proof}
The first sentence comes from a straightforward matrix rank calculation. In order to check the second sentence, we need to verify that $\pi_i\cdot V_i>0$ for $i\in\{1,2,3,4\}$ and $\pi_5\cdot \Omega>0$. The calculations are carried out in File~5 in~\cite{files}.
\end{proof}

According to Lemma~\ref{lempi14}, if $(\omega,v)$ is sufficiently close to $(\Omega,0,\ldots,0)$, then $\Delta(\omega,v)$ is still a simplex, and we refer to the facet of $\Delta(\omega,v)$ not containing $V_i(\omega,v)$ as the $i$th facet. Similarly, we choose $\pi_i(\omega,v)$ to be normalized vectors, and then $\pi_i(\omega,v)\cdot x\geqslant 0$ are the inequalities defining $\Delta(\omega,v)$. Also, we introduce the rational function $\Psi(\omega,v)$ defined by the formula $\pi_5(\omega,v)\cdot H$.

\begin{lem}\label{lempi5}
$\Psi$ is well-defined in some ball centered at $(\Omega,0,\ldots,0)$. We have
$$\Psi(\Omega,0,\ldots,0)=0,\,\,\,\,\,\,\,\,\left.\frac{\partial\Psi(\Omega(t),v)}{\partial t}\right|_{\substack{t=\tau \\ v=0}}=0,\,\,\,\,\,\,
\,\,\left.\frac{\partial^2\Psi(\Omega(t),v)}{\partial t^2}\right|_{\substack{t=\tau \\ v=0}}<0.$$
\end{lem}

\begin{proof}
As said above, our construction makes sense for all $(\omega,v)$ sufficiently close to $(\Omega,0,\ldots,0)$, which proves the first assertion. Further, we check that the function $\Psi(\Omega(t),0,\ldots,0)/\pi(t)$, where $\pi(t)$ is the polynomial as in Section~\ref{seccompdet}, is well-defined at $\tau$. Since $\pi$ has a double root at $t=\tau$, so does $\Psi(\Omega(t),0,\ldots,0)$, which proves the desired equalities. The inequality can be checked numerically. These computations are contained in File~4 of~\cite{files}.
\end{proof}

Further, we choose rational interior points $W,W_1\ldots,W_4\in\Delta$ satisfying $\|W-\Omega\|<\varepsilon$ and $\|W_i-V_i\|<\varepsilon$.
We define $P$ as the convex hull of $H$, $W$, $W_i$'s and $F_{ij}$'s. We also define $A_i,B_i,C_i,D_i$ as arbitrary rational points whose affine span is the hyperplane $\{x_i=0\}$ and whose convex hull has diameter less than $\varepsilon$ and contains $V_i$ as a relative interior point. Finally, we define $Q$ as the convex hull of $\Omega_1,\Omega_2$ and all $A_i$'s, $B_i$'s, $C_i$'s, $D_i$'s.

Now we can check the items~(I)--(III) in Theorem~\ref{prPolytope}. To check~(I), we note that the points $A_i,B_i,C_i,D_i,W_j$, $W$, $\Omega_1,\Omega_2$, $H$ are rational by their definitions, and the points $F_{ij}$ are rational functions in $\alpha,\beta,\gamma$. Since the vertices of $P$ and $Q$ are taken from these points, the item~(I) is valid.

We proceed with item (III) of Theorem~\ref{prPolytope}. The inclusion $\Delta\subset Q$ is straightforward because every vertex of $\Delta$ belongs to the convex hull of vertices of $Q$ directly by definitions. As to the inclusion $P\subset\Delta$, it is again straightforward to see that the points $W$, $W_i$ belong to $\Delta$. In order to prove $F_{ij}\in\Delta$, we check numerically that $\pi_k\cdot F_{ij}>0$ whenever $k\neq i$, and the fact that $\pi_i\cdot F_{ij}=0$, which means that $F_{ij}$ lies on the $i$th facet of $\Delta$, follows from our construction. Similarly, we check numerically that $\pi_k\cdot H>0$ if $k\neq 5$, and the equality $\pi_5\cdot H=0$ is a part of Lemma~\ref{lempi5}; so we have $H\in\Delta$. The computations of this paragraph are contained in File~5 in~\cite{files}.

In order to prove the item (II) of Theorem~\ref{prPolytope}, we note that the vertices $W, W_1, W_2, W_3, W_4$ of $P$ are close, respectively, to the points $\Omega,V_1,V_2,V_3,V_4$, whose convex hull is $\Delta$. This means that $\dim P=\dim\Delta=4$, and the fact that $P\subset Q$ is already established in the previous paragraph.

\section{Intermediate simplices}

This section is devoted to the item~(IV) in Theorem~\ref{prPolytope}. Namely, we are going to check that every simplex nested between $P$ and $Q$ does necessarily have a vertex not all of whose coordinates belong to $\F$. We have verified that $\tau\notin\F$ in Section~\ref{seccompdet}, so the simplex $\Delta$ satisfies this assumption.

Now let $\Delta'$ be a simplex satisfying $P\subset\Delta'\subset Q$. By taking a sufficiently small $\varepsilon$, we can guarantee that $\Delta'$ should be arbitrarily close to $\Delta$. In other words, we can assume without loss of generality that the vertices of $\Delta'$ belong to the balls centered at the vertices of $\Delta$ of radii not exceeding any positive number $\delta$ fixed in advance. We want to show that, if $\delta$ is sufficiently small, then $\Delta'=\Delta$, which would complete the proof of item~(IV) in Theorem~\ref{prPolytope}.

Since $\Delta'$ is close to $\Delta$, it has a vertex $\omega$ close to $\Omega$, and since $\Delta'\subset Q$ we should have $\omega\in Q$. Similarly, for any facet of $\Delta$, there is a facet of $\Delta'$ that is close to it. In particular, for any $i\in\{1,2,3,4\}$, the $i$th facet of $\Delta'$ passes through points $f_{i1}(v_{i1}), f_{i2}(v_{i2}), f_{i3}(v_{i3})$ with some small values of $v_{ij}$'s. In other words, the equation of the $i$th facet of $\Delta'$ is $\pi_i(\omega,v)\cdot x=0$, so the remaining vertices of $\Delta'$ lie on the straight line containing $\omega$ and $V_i(\omega,v)$. Since all the coordinates of $V_i$ are positive except that the $i$th coordinate is zero, the same is true for the coordinates of $A_i,B_i,C_i,D_i$, which means that $x_i\geqslant0$ is a facet defining inequality of $Q$. In other words, the $i$th vertex of $\Delta'$ lies in the convex hull of $\omega$ and $V_i(\omega,v)$, which implies $\Delta'\subset\Delta(\omega,v)$. 

As we can check (see File~5 in~\cite{files}), the point $f_{ij}(u_{ij})$ lies inside $\Delta$ if $u_{ij}$ is a small negative number. Therefore, since $\Delta(\omega,v)$ is close to $\Delta$, there is a $u_{ij}<0$ not depending on $\omega, v$ for which $f_{ij}(u_{ij})\in\Delta(\omega,v)$. Therefore, the inequality $v_{ij}<0$ would imply that the intersection of the line $f_{ij}(\R)$ with $\Delta(\omega,v)$ contains points $f_{ij}(w_{ij})$ with negative $w_{ij}$'s only, which contradicts to the fact that $f_{ij}(0)=F_{ij}\in P\subset\Delta'\subset\Delta(\omega,v)$. So we get $v_{ij}\geqslant 0$, and, since $H\in P\subset\Delta'\subset\Delta(\omega,v)$, we get $\Psi(\omega,v)\geqslant 0$.

Therefore, it remains to check that the inequality $\Psi(\omega,v)<0$ holds for all $(\omega,v)\in Q\times(\R_+)^{12}$ that are close to but not equal to $\Omega_0=(\Omega,0,\ldots,0)$. 

\begin{obs}\label{lemmax}
Let $\Phi(t,z_1,\ldots,z_n)$ be a rational function on $\R\times\R^n$ that takes a zero value at the point $O=(0,\ldots,0)$. Assume that $\partial\Phi/\partial z_i(O)<0$ holds for all $i$, and also $\partial\Phi/\partial t(O)=0$, $\partial^2\Phi/\partial t^2(O)<0$. Then $\Phi(\xi,\zeta_1,\ldots,\zeta_n)<0$ holds for all $(\xi,\zeta_1,\ldots,\zeta_n)$ that are different from but sufficiently close to $O$ and satisfy $\zeta_i\geqslant0$.
\end{obs}

\begin{proof}
Follows from Taylor's theorem.
\end{proof}

The function $\Psi$ is rational and satisfies $\Psi(\Omega_0)=0$. Denoting by $l_1,\ldots,l_{16}$ the directions of the rays pointing from $\Omega$ to $(A_i, B_i, C_i, D_i)$, we set
$$\Phi(t,z_1,\ldots,z_{16},v)=\Psi(\Omega(t)+z_1l_1+\ldots+z_{16}l_{16},v).$$
In order to complete our proof, we need to check that the assumptions of Observation~\ref{lemmax} hold for $\Phi$. Namely, we check in File~5 of~\cite{files} that $\partial\Psi/\partial v_{ij}(\Omega_0)<0$, we confirm in File~6 that the derivatives of $\Psi$ along $\overrightarrow{\Omega V_i}$ are negative at $\Omega_0$, and since $A_i, B_i, C_i, D_i$ are close to $V_i$, this implies $\partial\Phi/\partial z_k(\Omega_0)<0$. Finally, we conclude from Lemma~\ref{lempi5} that $\partial\Phi/\partial t(\Omega_0)=0$, $\partial^2\Phi/\partial t^2(\Omega_0)<0$.

\section{Conclusion}

We have constructed a subfield $\F\subset\R$ and a pair of polytopes $P\subset Q\subset\R^{d-1}$ with vertices in $\F^{d-1}$ such that the existence of a $(d-1)$-simplex $\Delta$ nested between $P$ and $Q$ depends on the extension of $\F$ which is allowed to contain the vertices of $\Delta$. In terms of linear algebra, this gives a matrix $A$ with nonnegative entries in $\F$ satisfying
\begin{equation}\label{eq6}
d=\mathrm{rank}(A)=\mathrm{rank}_+(A,\R)<\mathrm{rank}_+(A,\F).
\end{equation}
In our construction, we have $d=5$, and by constructing block-diagonal matrices we can get the same result with any $d>5$. The condition~\eqref{eq6} is impossible for $d\leqslant 3$ (see~\cite{KRS}), but we do not know what happens when $d=4$.

\begin{prob} (See also~\cite{anotherproof, KK, Consp}.)
Does there exist a real $m\times n$ matrix $A$ such that $\mathrm{rank}_+(A,\R)=4$ and $\mathrm{rank}_+(A,\mathbb{Q}(a_{11},\ldots,a_{mn}))=5$? 
\end{prob}

In our construction, $\mathcal{F}$ is a finite extension of $\mathbb{Q}$, and it would be interesting to know if the condition~\eqref{eq6} is possible with $\F=\mathbb{Q}$. We suspect that stronger results, similar to the universality theorem for nonnegative factorizations~\cite{myuni}, may hold for \textsc{exact nmf} as well. In particular, we expect that matrices satisfying~\eqref{eq6} exist over any field $\F$ that is not real closed, and we conjecture that \textsc{exact nmf} is an $\exists\R$-complete problem. Computational difficulties arising already in the simplest cases as the one considered in the present paper make us afraid, however, that such a conjecture does not admit a proof of reasonable length. A particular problem that seems to be beyond our reach is to construct a rational matrix $A$ satisfying
$$\mathrm{rank}(A)=\mathrm{rank}_+(A,\mathbb{Q})<\mathrm{rank}_+(A,\mathbb{Q}(\sqrt[p]{2})),$$
where $p$ is a sufficiently large prime (e.g., $p=3$).

\section{Acknowledgments}

I am grateful to Kaie Kubjas for an interesting discussion on the topic and helpful suggestions on the presentation of the result. A part of this discussion held in 2015 in \textit{Aalto University} in Helsinki, and I am grateful to Kaie for inviting me and to the colleagues from the university for their hospitality. I would like to thank Till Miltzow for a discussion and for pointing me to~\cite{Consp}. I am grateful to the Editor-in-Chief and Associate Editor who were handling my submission for very detailed specific comments concerning its content and to anonymous reviewers for careful reading, valuable comments, and kind words of encouragement.

\end{document}